# An aperiodic monotile for the tiler

Vincent Van Dongen

## Abstract
*Can the entire plane be paved with a single tile that forces aperiodicity? This is known as the ein Stein problem (in German, ein Stein means one tile). This paper presents an aperiodic monotile for the tiler. It is based on the monotile developed by Taylor and Socolar (whose aperiodicity is forced by means of a non-connected tile that is mainly hexagonal) and motif-based hexagonal tilings that followed this major discovery. The proposed monotile consists of two layers. No motif is needed to make the monotile aperiodic. Additional motifs can be added to the monotile to provide some insights. The proof of aperiodicity is presented with the use of such motifs.*

## Introduction

Penrose discovered a set of two tiles with some amazing properties (Penrose, 1979). Each tile of this set cannot pave alone the entire plane. And when used together, the two tiles can only generate aperiodic tilings. Soon after, a fundamental question was raised. Could one tile alone be aperiodic? That is, could a monotile pave the entire plane while not allowing any periodic tiling? This problem is known as the *'ein-Stein'* problem. ('ein Stein' is German and means 'one tile'.)

About a decade ago, a first solution to this problem was proposed (Socolar & Taylor, 2011). The tile is mainly hexagonal but with non-connected borders to enforce the aperiodicity. The tile and its flip side are shown here below as well as the tiling it generates.

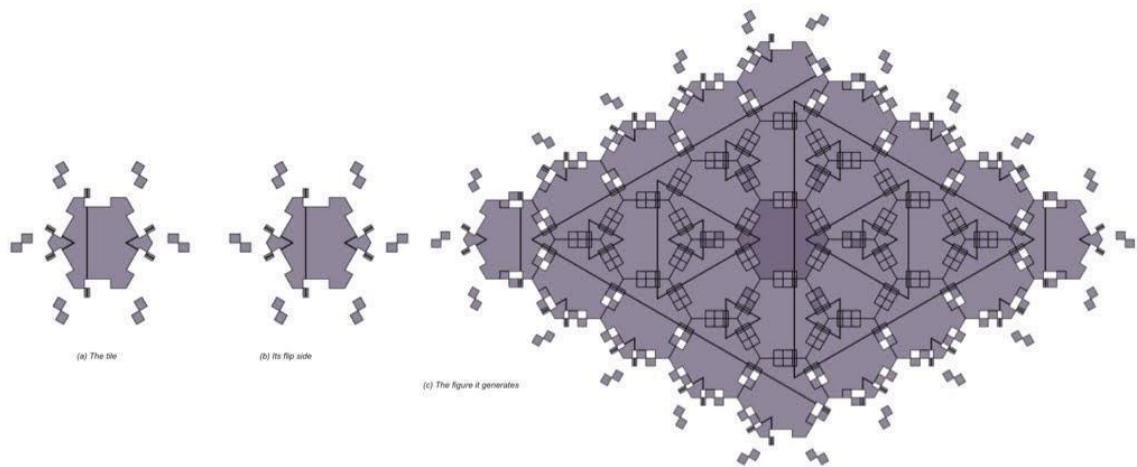

*Figure 1: Non-connected tile proposed as a solution to the ein-Stein problem and the Sierpinski triangle it generates when using it.*

Note that the motif on the tile (the black lines) is not necessary as the shape of the tile alone creates the aperiodicity. There is only one way of covering the plane with it and it requires both sides of the tile (the other side is obtained by flipping the tile). The black lines on the tile provide some insight on the tiling being generated. On the paved area, a figure gets created that is a fractal similar to the well-known triangle of Sierpinski (Sierpinski triangle). See figure here above on the right.

This tile is unfortunately not connected. In other words, it is not in one piece. This makes it impossible for the tiler to use it in practice. Other monotiles were proposed with a motif on the tile that enforces aperiodic tiling. For example, a motif on hexagonal tiles was recently proposed in (Mampusti & Whittaker, 2020). The solution makes use of a dendrite motif. As explained in their paper, the dendrite forms nested triangles of Sierpinski as well. Their innovation is that the monotile does not require to be flipped and the enforcement of the aperiodic property of the tiling is controlled by the dendrite motif. In a dendrite, there is no cycle and this feature allows them to prove the aperiodic character of the tiling.

Their proposal was the basis of a self-ruling monotile (Gradit & Van Dongen, 2022) that forces aperiodic tiling by design; the way of using the monotile is defined in the tile itself. However, the shape of the tile alone doesn't make the tile aperiodic.

In this paper, we present yet another version of a hexagonal-based tile called **HexToo**. It is aperiodic in the sense that it cannot be used to generate periodic tiling. Like in TS-tiling, it is only the shape of this new monotile that enforces aperiodic tiling.

Like in the dendrite, our monotile is such that no periodic cycle can exist. A simple trick was applied to its shape to avoid such cycle: **HexToo** consists of two layers with a bottom part that can only be covered by a top part. We show how a tiler can use **HexToo** to pave the plane: tile placement ordering is key to produce a successful tiling. Hence, a well-educated tiler will be able to tile any surface with **HexToo**, and the result is aperiodic.

## A two-layer monotile called HexToo

Our monotile called **HexToo** is shown here below. It is made of two layers: (i) a top one that is the visible part of the tile, and a bottom (ii) that it is underneath.

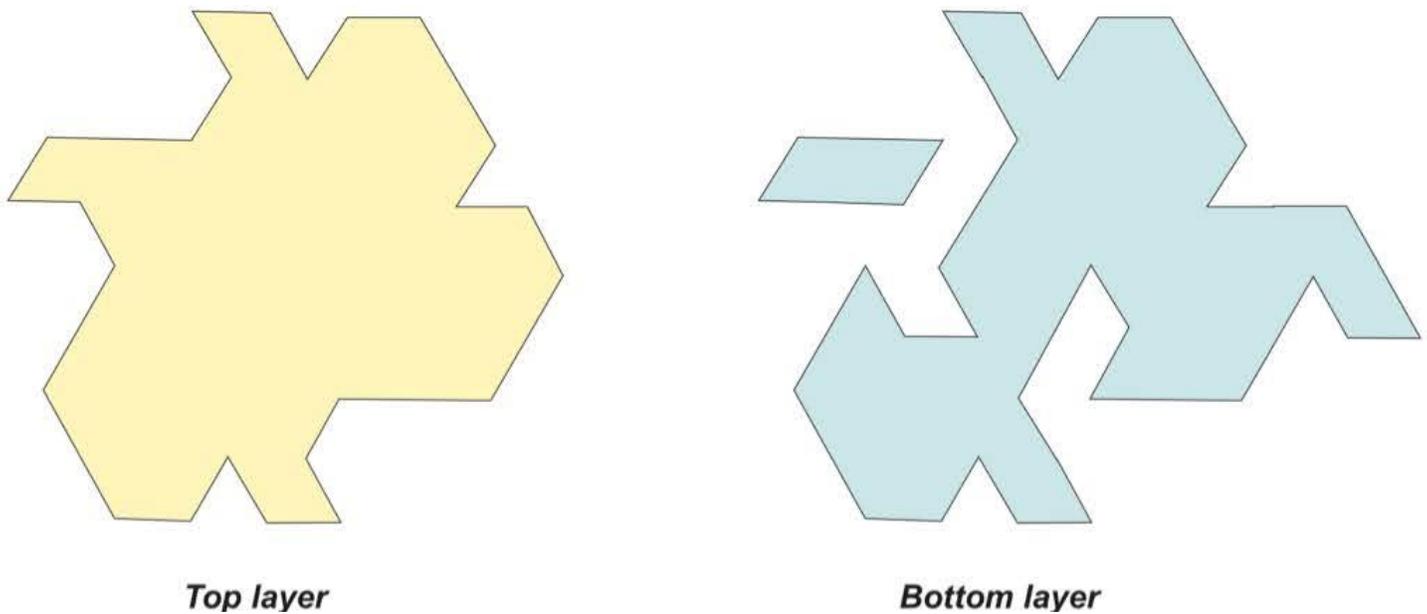

*Figure 1: HexToo is an aperiodic monotile made of two thin layers.*

The following figure shows both sides of the monotile: recto and verso (or flip side). The two layers are there to create joints. The monotile offers one male joint and available female joints. The female joints are not always used but are always available.

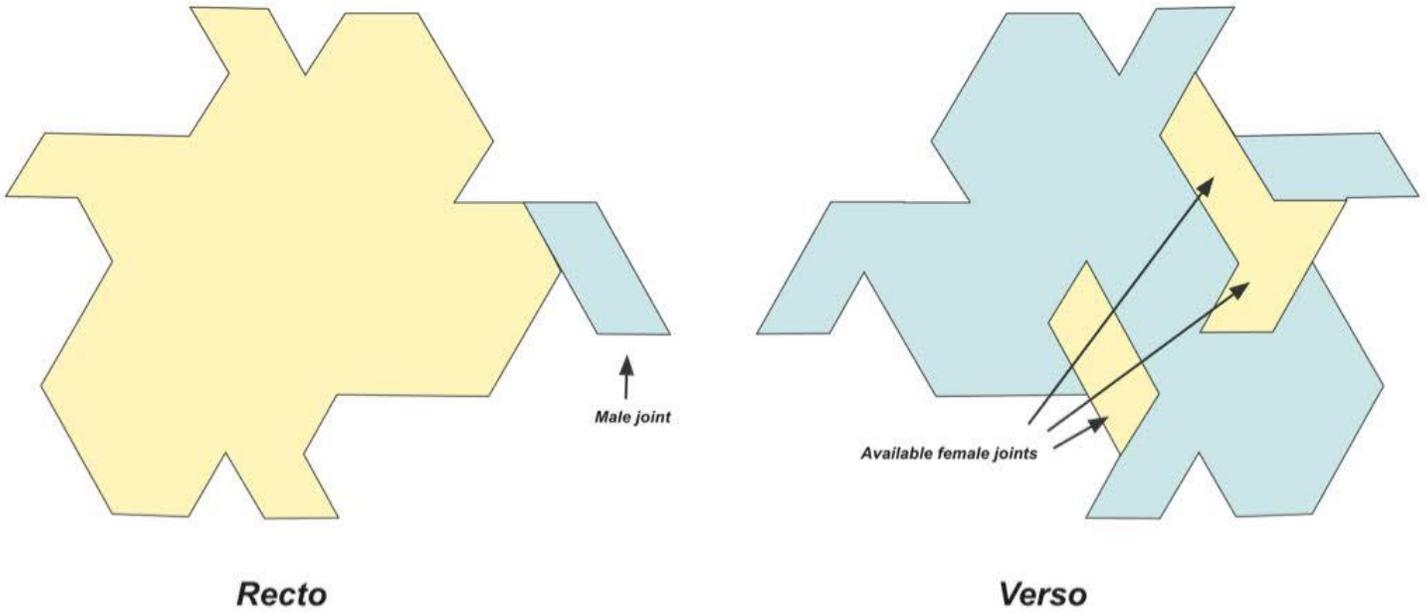

*Figure 2: Both sides of the monotile: recto and verso (flip side). The two layers define male and female joints.*

The following figure shows a small tiling made of 19 tiles. The flip side of the tiling shows how tiles fit into each other, thanks to the joints.

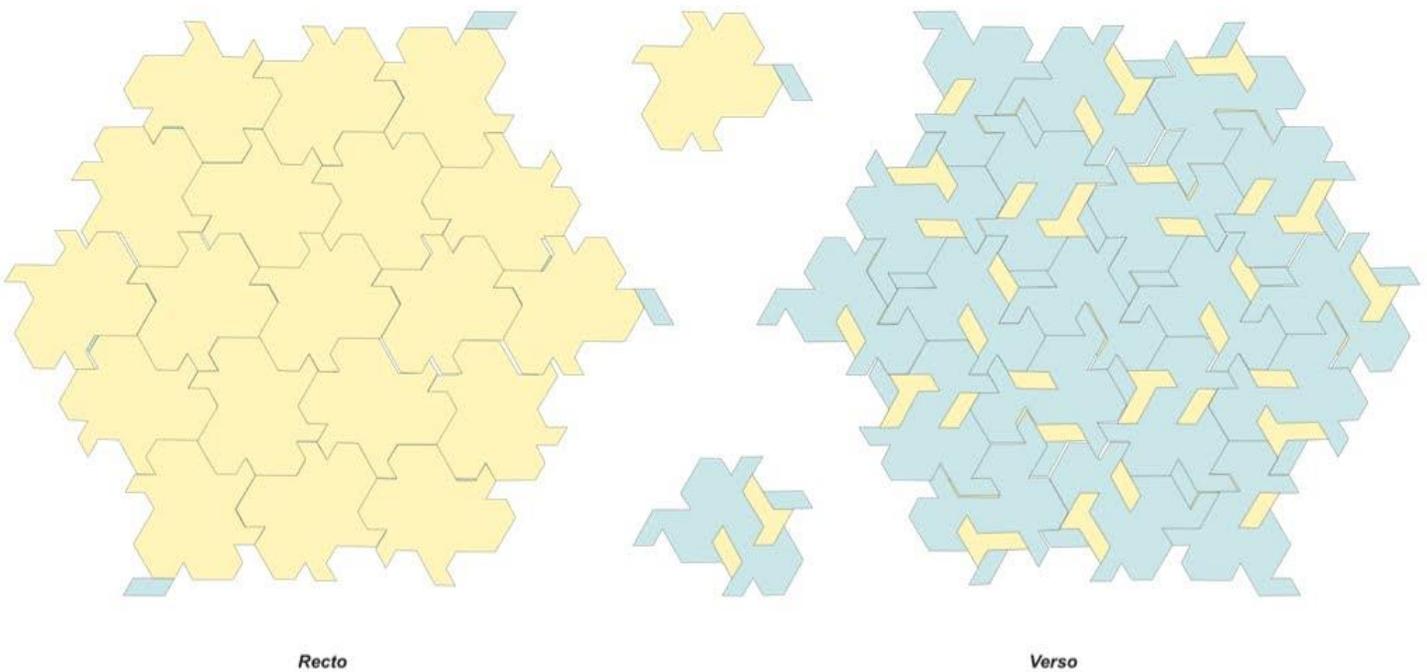

*Figure 3: A tiling made of HexToo. The flip side shows how tiles fit into each other, thanks to the joints.*

The shape of **HexToo** enforces the tiling to be aperiodic. The best way to see this is by adding a motif on the tile. See the figure here below.

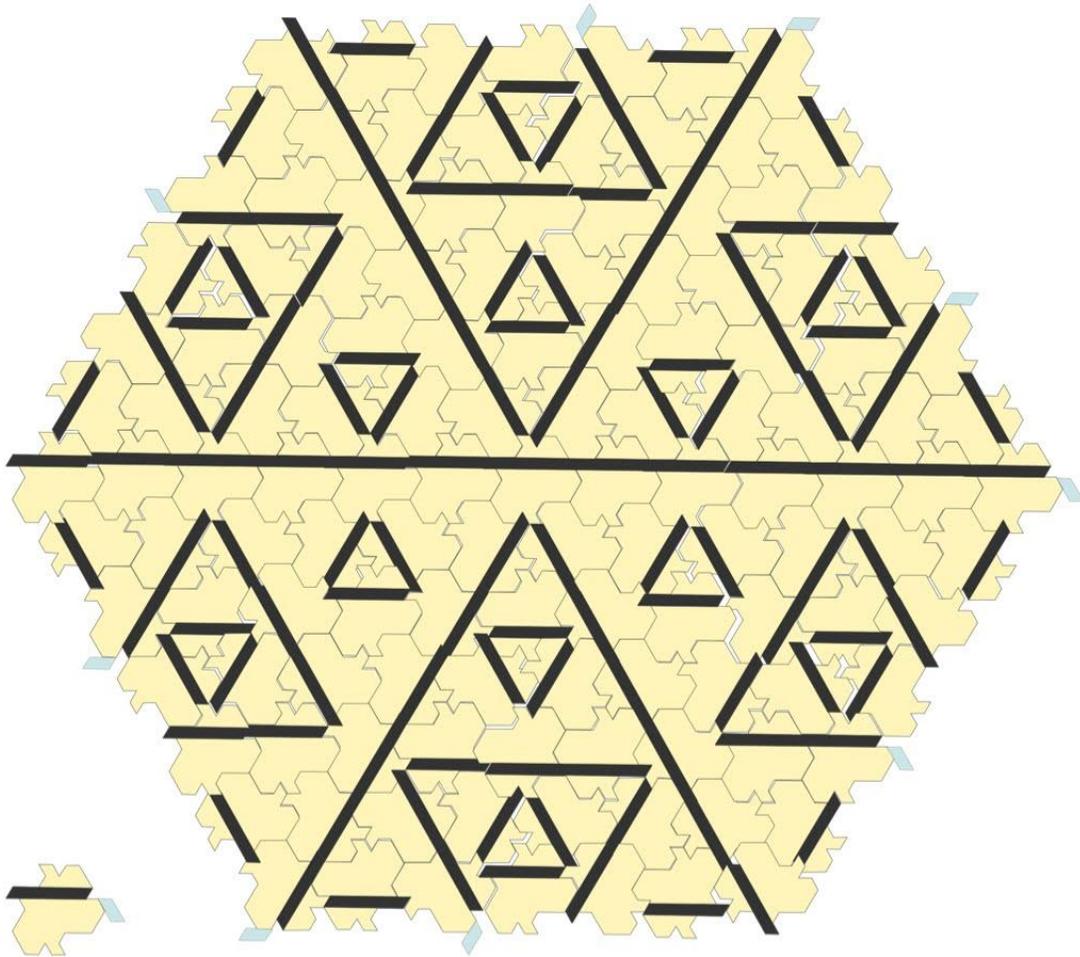

*Figure 4: A motif on HexToo shows its aperiodic property, Sierpinski triangles being created.*

The motif created are Sierpinski triangles like in previously discovered tilings, as in (Socolar & Taylor, 2011), (Mampusti & Whittaker, 2020) and (Gradit & Van Dongen, 2022).

We still need to demonstrate that **HexToo** is aperiodic, i.e. only aperiodic tilings can be created. For doing so, we will make use of another motif: the dendrite motif introduced in (Mampusti & Whittaker, 2020).

## About the aperiodic character of the monotile HexToo

In (Mampusti & Whittaker, 2020), a dendrite motif is used to force the aperiodicity of the tiling. The demonstration is that the tiling is aperiodic as long as no cycle is created (by definition, there is no cycle in a dendrite).

In this section, we present a motif on **HexToo** for the dendrite. We then show why **HexToo** will never create a cycle in the dendrite. Hence, this will prove the aperiodic character of **HexToo**.

**HexToo** with the dendrite motif is shown here below along with a tiling.

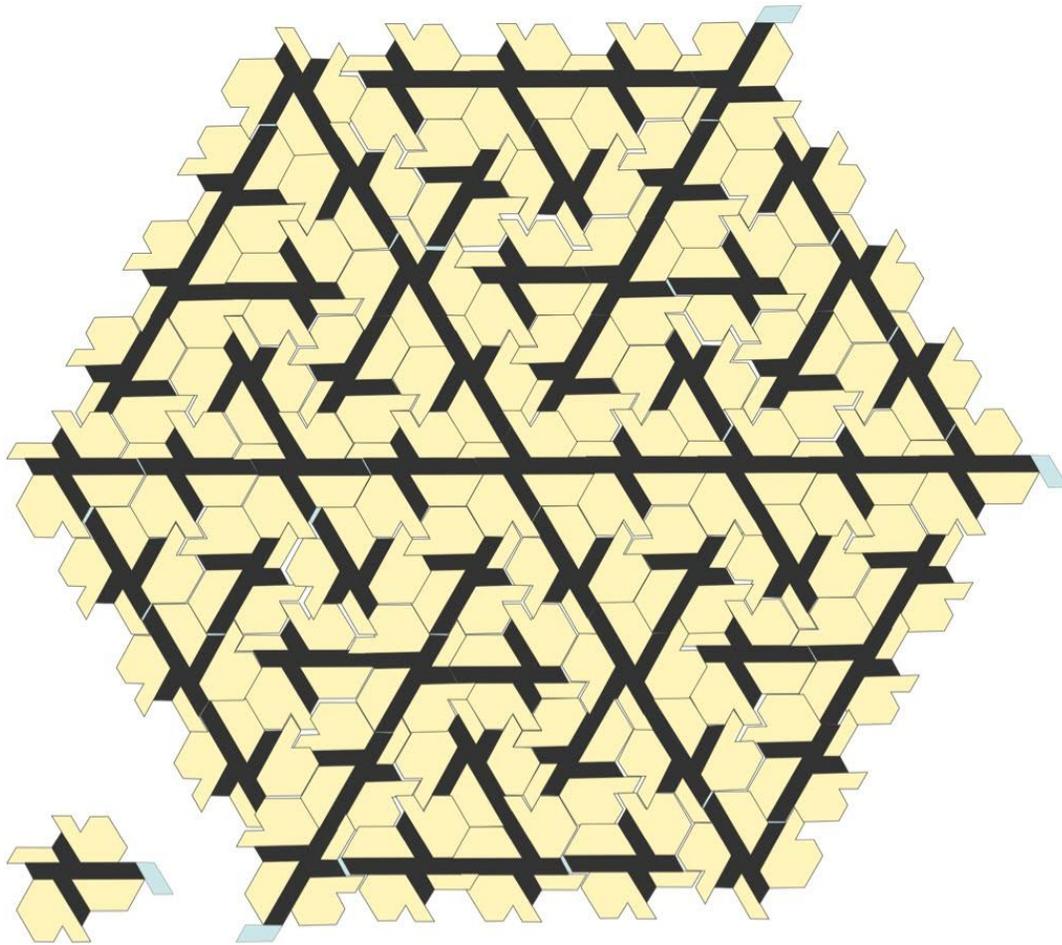

*Figure 5: HexToo with a motif that produces a tiling with dendrites (i.e. without cycle).*

In (Mampusti & Whittaker, 2020), the motif on the monotile enforces the aperiodicity of the tiling. With **HexToo**, the shape of the monotile, with its two layers and joints, such a cycle is not possible, as shown here below.

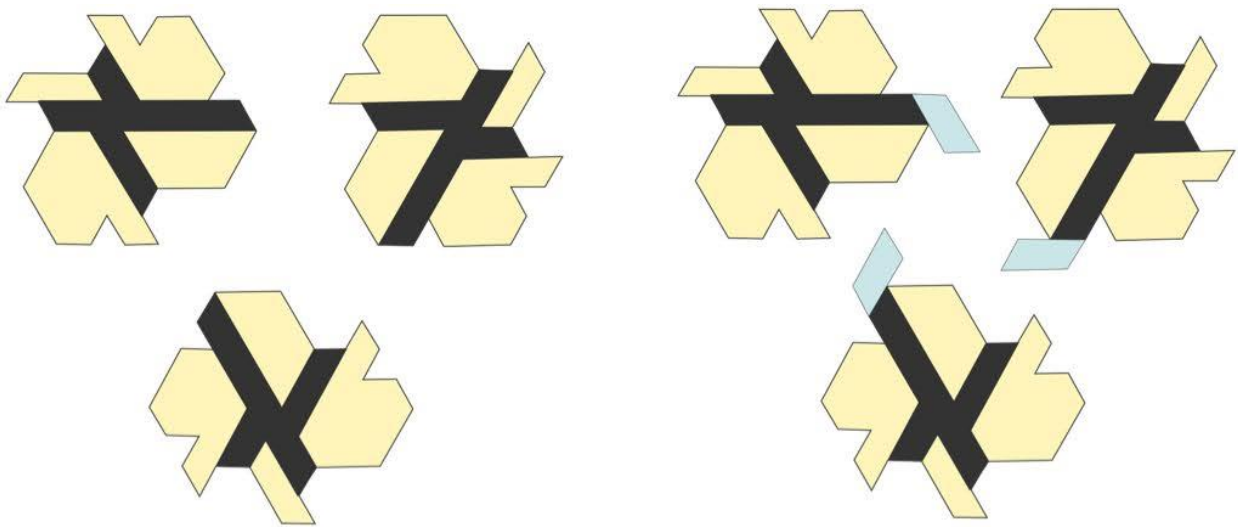

*Figure 6: The two-layers of HexToo enforce the aperiodic tiling as dendrites will have no cycle by design.*

The no-cycle constraint for aperiodicity works because (i) this motif on HexToo is equivalent to the one on the monotile of (Mampusti & Whittaker, 2020) shown here below on the left, and (ii) both tiles are hexagonal with inner and outer bumps at the same edges to implement polarity.

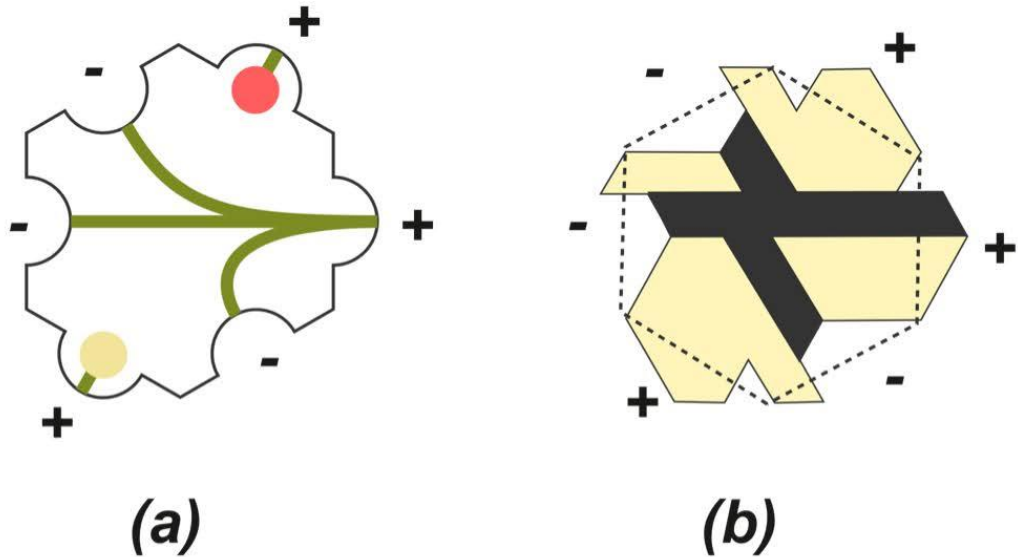

*Figure 7: Two monotiles with similar shapes and motifs.*

## How to use HexToo (i.e. instructions for the tiler)

The two layers of **Hextoo** force tiles to be placed in a specific order. As explained in the previous section, no order exists for a cycle to be created and this makes the monotile aperiodic. The tiler needs to be careful though. Because of its two layers, tiles must be placed in the right sequence to properly cover the area. (Tiles that consist of one layer can be placed in any order.)

The good news for the tiler is that an order does exist for the tiles to pave any given area. Tiles can be placed one tile at a time. The sequence must follow the motif of the dendrite, starting from its endpoints. The following figure shows how a motif on the tile can create a dependency graph for tile ordering.

.

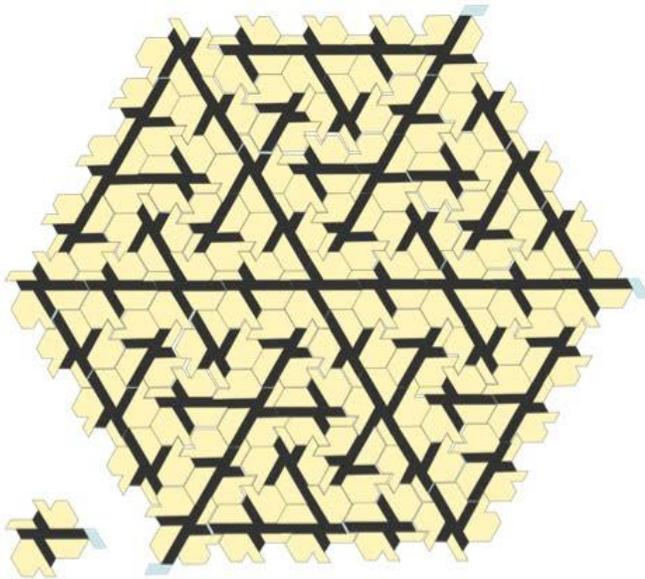
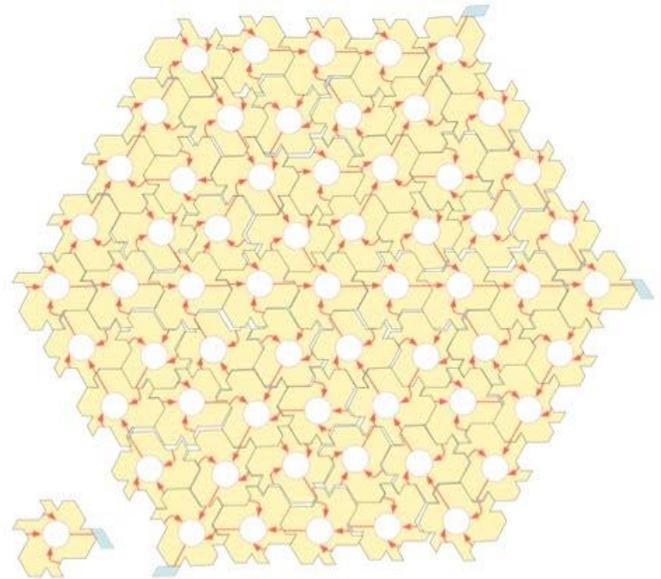

(a) Dendrite motif.

(b) Motif for dendrite to become a dependency graph for tile ordering.

*Figure 8: Dendrite motif and special motif for creating a dependency graph for tile ordering..*

## Alternate motifs

As explained above, it is the shape alone of the monotile **HexToo** that makes it aperiodic and this shape consists of two layers. Any motif can be printed on **HexToo**. Here below is an example of yet another dendrite motif and its resulting tiling.

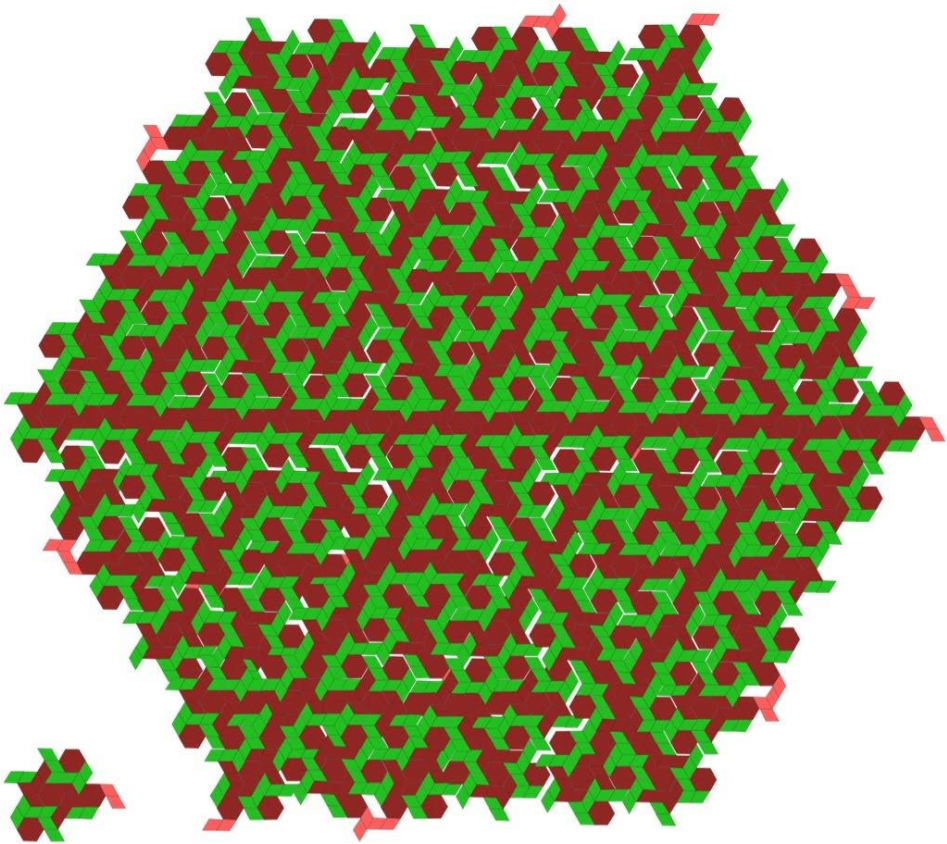

*Figure 9: HexToo with yet another dendrite motif on it and the corresponding aperiodic tiling.*

Here below is yet another motif where the grey areas mark the position of the female joints on the lower layer.

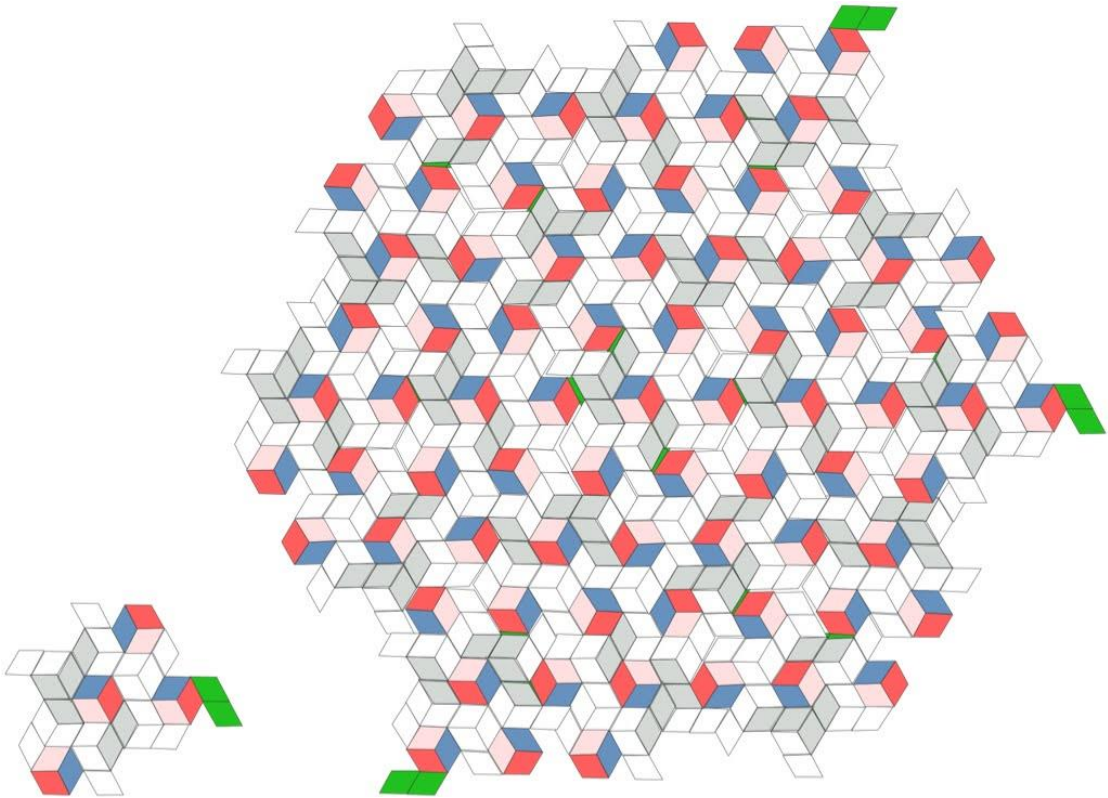

*Figure 10: HexToo with a motif using the same rhombus but in different colors, and the corresponding tiling.*

## Alternate shapes

**HexToo** is a two-layer monotile that is hexagonal-base. Different variations of such a monotile can be created. In fact, there is an infinite number of them. See here below a few of them to illustrate our point; the area changes but remains hexagonal-base. Please note how some rhombuses change shape between (a), (b) and (c), while others remain the same.

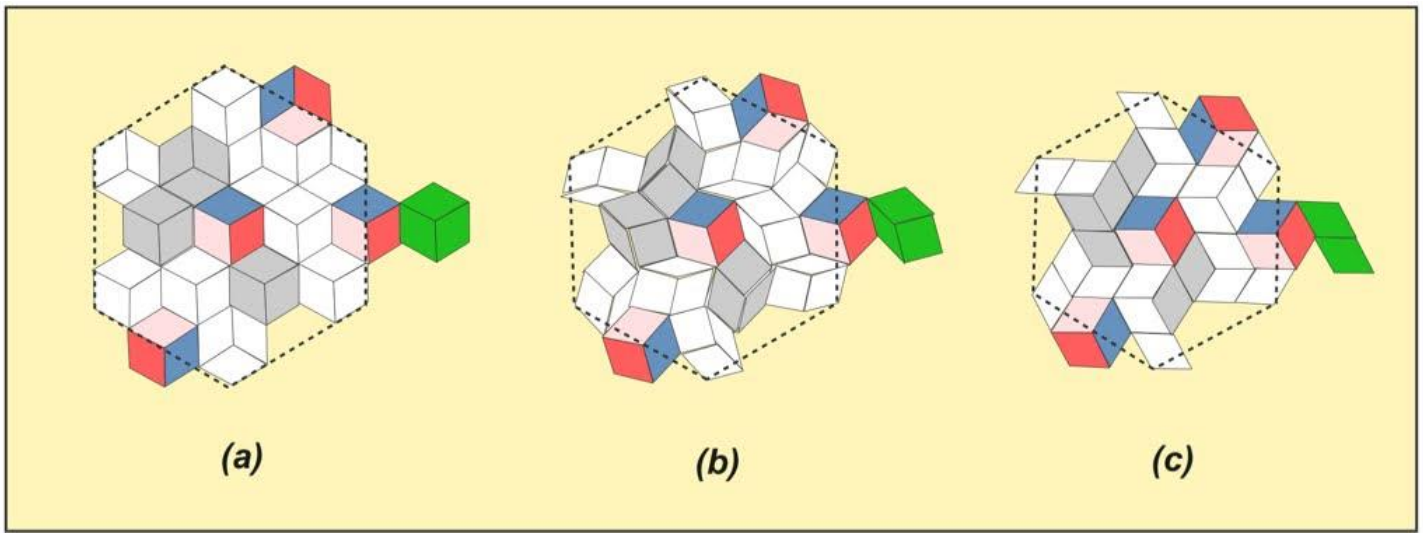

# Conclusion

In this paper, we presented a new monotile called **HexToo** whose shape, and only its shape, makes it aperiodic. It is hexagonal-base and consists of two layers with male and female joints. The reason why it is aperiodic is provided. Its proof makes use of the result presented in (Mampusti & Whittaker, 2020): dendrite motifs (with no cycle) are being enforced by the layers and joints of **HexToo**.

With these layers, tiles must be placed in the right order. We show that a correct sequence always exist. A tiler can follow such a sequence to pave any given area.

Like the Socolar-Taylor tile, a motif on **HexToo** can create a tiling with Sierpinski triangles.

At the end, some other motifs for **HexToo** are presented as well as valid changes to its hexagonal-base shape.

# References

*Einstein Problem*. (n.d.). Retrieved from Wikipedia: https://en.wikipedia.org/wiki/Einstein_problem

Mampusti, M., & Whittaker, M. F. (2020). An aperiodic monotile that forces nonperiodicity through dendrites. University of Wollongong, Australia.

*Sierpinski triangle.* (n.d.). Retrieved from Wikipedia: https://en.wikipedia.org/wiki/Sierpi%C5%84ski_triangle

Gradit, P., & Van Dongen, V. (2022). *A self-ruling monotile for aperiodic tiling.* Arxiv: https://arxiv.org/ftp/arxiv/papers/2201/2201.03079.pdf

Penrose, R. (1979). Pentaplexity A Class of Non-Periodic Tilings of the Plane. *The Mathematical Intelligencer 2*, 2–37.

Socolar, J., & Taylor, J. (2011). An aperiodic hexagonal tile. *Journal of Combinatorial Theory, Series A,*, 2207–2231.

*Socolar-Taylor tile.* (n.d.) Retrieved from Wikipedia: https://en.wikipedia.org/wiki/Socolar%E2%80%93Taylor_tile